\documentclass[11pt]{article}

\usepackage[centering,margin=1.2in]{geometry}

\usepackage{datetime}

\usepackage{tikz}

\usepackage{eucal}

\usepackage[all,cmtip]{xy}

\usepackage{amssymb, amsmath, amsthm}

\usepackage[shortlabels]{enumitem}

\numberwithin{equation}{section}

\newtheorem{theorem}{Theorem}[section]
\newtheorem{proposition}[theorem]{Proposition}

\newtheorem{lemma}[theorem]{Lemma}

\newtheorem*{theorem-without-number}{Theorem}

\theoremstyle{definition}

\newtheorem{definition}[theorem]{Definition}

\newtheorem{example}[theorem]{Example}
\newtheorem{remark}[theorem]{Remark}

\usepackage{todonotes}

\usepackage[pdftex]{hyperref}

\hypersetup{
  colorlinks = true,
  linkcolor = black,
  urlcolor = gray,
  citecolor = black,
}

\newcommand{\rmB}{\mathrm{B}}
\newcommand{\rmH}{\mathrm{H}}
\newcommand{\rmK}{\mathrm{K}}
\newcommand{\rmS}{\mathrm{S}}

\newcommand{\bbA}{\mathbb{A}}
\newcommand{\bbC}{\mathbb{C}}
\newcommand{\bbP}{\mathbb{P}}
\newcommand{\bbS}{\mathbb{S}}
\newcommand{\bbQ}{\mathbb{Q}}
\newcommand{\bbZ}{\mathbb{Z}}

\newcommand{\calB}{\mathcal{B}}
\newcommand{\calG}{\mathcal{G}}
\newcommand{\calR}{\mathcal{R}}

\newcommand{\Pic}{\mathrm{Pic}}
\newcommand{\Bir}{\mathrm{Bir}}

\newcommand{\Burn}{\mathrm{Burn}}
\newcommand{\Aut}{\mathrm{Aut}}
\newcommand{\Cr}{\mathrm{Cr}}
\newcommand{\PGL}{\mathrm{PGL}}
\newcommand{\GL}{\mathrm{GL}}
\newcommand{\Hom}{\mathrm{Hom}}

\newcommand{\id}{\mathrm{id}}
\newcommand{\sign}{\mathrm{sign}}
\newcommand{\pr}{\mathrm{pr}}

\newcommand{\GGa}{\mathbb{G}_\mathrm{a}}
\newcommand{\GGm}{\mathbb{G}_\mathrm{m}}

\DeclareMathOperator{\colim}{colim}

\renewcommand{\phi}{\varphi}

\title{\bf Homological stability fails for the Cremona groups, yet their stable homology is non-trivial}

\author{Markus Szymik}

\newdateformat{mydate}{\monthname~\twodigit{\THEYEAR}}
\date{\mydate\today}

\begin{document}

\maketitle


The Cremona groups are the groups of all birational transformations of rational varieties, or, in other words, the groups of all automorphisms of rational function fields. These groups are traditionally studied one dimension at a time. Instead, we will consider the whole sequence as the dimension increases. This shift in perspective raises questions about `convergence' and the~`limit'. Contrary to the standard expectation for families of automorphism groups, we show that the induced homomorphisms in homology are rarely injective or surjective. This means that homological stability does not hold for the Cremona groups. We can nevertheless consider their stable homology, and we show it is represented by a non-trivial infinite loop space.






\section{Introduction}

The Cremona groups are groups of coordinate transformations that play a central role in algebra and geometry. Algebraically, for any integer~$n\geqslant0$ and any field~$F$, the Cremona group~$\Cr_n(F)$ is the group~$\Aut(F(x_1,\dots,x_n)|F)$ of~$F$--automorphisms of the rational function field~$F(x_1,\dots,x_n)$ in~$n$ variables~\hbox{$x_1,\dots,x_n$}. Geometrically, if~$X$ is an~$n$--dimensional rational~$F$--variety such as projective space~$\bbP^n_F$, 
affine space~$\bbA^n_F=(\GGa^n)_F$, or a torus~$(\GGm^n)_F$, we can identify the group~$\Bir(X)$ of all birational self-equivalences of~$X$ with~$\Cr_n(F)$ via the induced action on its function field~\hbox{$F(X)\cong F(x_1,\dots,x_n)$}, and the identification is well-defined up to conjugation. The Cremona groups have been the subject of intensive studies for over a century. 

Nevertheless, it seems fair to say that they are far from well understood beyond some low-dimensional, exceptional cases, and that very little was known about the Cremona groups in higher dimensions. However, recent progress suggests the time is ripe to change this~(see~\cite{BLZ, BSY, GLU, Lin+Shinder:1} and the references therein). 

In this paper, we leave the traditional perspective, which considers each dimension at a time, behind and instead ask whether any systematic patterns might appear as we let~$n$ tend to infinity. More precisely, given an element~$\phi\colon x=(x_1,\dots,x_n)\mapsto\phi(x)$ in~$\Cr_n(F)$, we can add a new variable~$y$ and produce a new element~\hbox{$(x,y)\mapsto(\phi(x),y)$} in the next group~$\Cr_{n+1}(F)$. This gives rise to a diagram
\begin{equation}\label{eq:Cremona_diagram}
\dots\longrightarrow\Cr_n(F)\longrightarrow\Cr_{n+1}(F)\longrightarrow\dots
\end{equation}
of groups and morphisms between them. Of course, we could also have placed the new variable at the front or in the middle, and the resulting group morphisms would have changed by conjugation.

A standard question to ask for a diagram of groups such as~\eqref{eq:Cremona_diagram} is whether it satisfies homological stability~(see~\cite{RWW} and~\cite{Wahl} for a recent survey). Homology is a method of studying groups using linear algebra as the derived functors of abelianisation. It originated in the early works on class field theory, and the description of the Cremona groups as automorphism groups of field extensions indicates that we remain close in spirit to the Galois cohomology context. A group's first and second homology contain valuable information about generators and relations, respectively. Families of groups that satisfy homological stability abound, such as the symmetric groups~\cite{Nakaoka}, the braid groups~\cite{Arnold}, other mapping class groups~\cite{Harer}, and the general linear groups over fields~\cite{Quillen:1}. Aspects of homological stability related to function fields also feature in~\cite{EVW}. 

Given the ample evidence, it seems reasonable to expect that many more families of automorphism groups satisfy homological stability. Be the psychology of expectations as it may, we should be able to decide whether or not this property holds for such a prominent family as the Cremona groups. This paper presents a proof that homological stability fails for them. To summarise Theorems~\ref{thm:not_injective} and~\ref{thm:not_surjective} in the main text:

\begin{theorem}
For all~$n\geqslant 3$, the stabilisation homomorphism
\[
\rmH_1(\Cr_n(\bbC);\bbZ)\longrightarrow\rmH_1(\Cr_{n+1}(\bbC);\bbZ)
\]
has a non-trivial kernel and cokernel.
\end{theorem}

We start in Section~\ref{sec:homological_instability} with some elementary and natural examples that illustrate the general phenomenon. In Section~\ref{sec:inj_fails}, we show that the homomorphism induced by the group morphism~$\Cr_n(\bbC)\to\Cr_{n+1}(\bbC)$ in homology need not be injective, and in Section~\ref{sec:sur_fails}, we show that it need not be surjective. Nevertheless, we can still consider the stable homology, as it is a formal colimit, and discuss it in Section~\ref{sec:stable}. We summarise Theorems~\ref{thm:stable} and~\ref{thm:non-trivial}:

\begin{theorem}
For any field~$F$, the stable homology~$\colim\limits_n\rmH_*(\Cr_n(F);\bbZ)$ of the Cremona groups~{\upshape\eqref{eq:Cremona_diagram}} is given as the homology of an infinite loop space. If the field is infinite, that space is not discrete, and the stable homology is non-zero even rationally.
\end{theorem}

Throughout this text, we treat the Cremona groups as discrete groups. For locally compact fields, there is an interesting way to equip them with a Hausdorff topology that makes them connected topological groups; see~\cite{BF13, BZ18}. However, this is not our concern here.


\section*{Acknowledgements} 

The author thanks Serge Cantat, Julia Schneider, Evgeny Shinder, Burt Totaro, and an anonymous referee for their generous discussions that led to the account given here. The author would also like to thank the Isaac Newton Institute for Mathematical Sciences, Cambridge, for support and hospitality during the programme~`Equivariant Homotopy Theory in Context', where work on this paper was undertaken. This work was supported by EPSRC grant EP/Z000580/1.


\section{Homological instability}\label{sec:homological_instability}

We begin by recalling the definition.

\begin{definition}
A diagram
\begin{equation}\label{eq:diagram}
G_0\longrightarrow G_1\longrightarrow G_2\longrightarrow\dots
\end{equation}
of groups satisfies~\textit{homologically stability} if, for each degree~$d$, the induced diagram
\[
\rmH_d(G_0)\longrightarrow\rmH_d(G_1)\longrightarrow\rmH_d(G_2)\longrightarrow\dots
\]
consists almost always of isomorphisms.
\end{definition}

If a diagram is homologically stable, we would also like to know the range of morphisms that induce isomorphisms in homology, but since we are interested in diagrams that are not homologically stable, this is not relevant here.

The symmetric groups~$\rmS_n$~\cite{Nakaoka}, the general linear groups over fields~\cite{Quillen:1} and nice rings, the braid groups~\cite{Arnold}, the Higman--Thompson groups~\cite{Szymik+Wahl}, the automorphism groups of free groups~\cite{Hatcher}, 
free nilpotent groups~\cite{Szymik}, and other mapping class groups~\cite{Harer} all satisfy homological stability for a canonical choice of stabilisation maps between them: adding a generator that is fixed by the smaller group. In fact, the canonical choice is only well-defined up to conjugation, but conjugate group morphisms induce homotopic maps between the canonical resolutions, so that they induce the same maps in homology. 

In contrast, whenever the homology of a group~$G$ is non-trivial, the homology is unstable for the diagram of trivial morphisms between the groups~$G_n=G$. We see that homological stability is a property of the entire diagram, not just the family of groups. The examples below~(see Propositions~\ref{prop:Boole} and~\ref{prop:GL}) illustrate this phenomenon well. They are based on the following lemma.

\begin{lemma}\label{lem:det}
Let~$f\colon V\to V$ and~$g\colon W\to W$ be endomorphisms of finite-dimensional vector spaces~$V$ and~$W$. Then
\[
\det(f \otimes g)=\det(f)^{\dim(W)}\det(g)^{\dim(V)}.
\]
\end{lemma}

\begin{proof}
By symmetry and
\[
\det(f \otimes g)=\det(f \otimes\id_W)\det(\id_V \otimes g),
\]
we can assume that~$g=\id_V$. Then, we can use the relation
\[
\det(f\otimes\id_W)=\det(f)^{\dim(W)},
\]
which is clear, as~$f\otimes\id_W$ has a matrix representation of~$\dim(W)$ many blocks consiting of a copy of~$f$ each.
\end{proof}

\begin{proposition}\label{prop:Boole}
Let 
\[
\rmS_{2^0}\longrightarrow\rmS_{2^1}\longrightarrow\rmS_{2^2}\longrightarrow\dots
\]
be the diagram of the symmetric groups, where the morphisms are induced by taking a permutation~$\pi$ to~$\pi\times\id_2$, where~$\id_2$ is the identity permutation on a two-element set. Then the induced stabilisation homomorphisms
\[
\rmH_1(\rmS_{2^n})\to\rmH_1(\rmS_{2^{n+1}})
\]
are never surjective and not injective for~$n\geqslant 1$.
\end{proposition}

\begin{proof}
The sign of a permutation is the determinant of the corresponding permutation matrix. Therefore, Lemma~\ref{lem:det} gives
\[
\sign(\pi\times\id_2)=\sign(\pi)^2=1,
\]
regardless of the permutation~$\pi$. This shows that we have commutative diagrams
\[
\xymatrix@C=4em{
\rmS_{2^n}\ar[r]^-{(-)\times\id_2}\ar[d]_-{\text{sign}} & \rmS_{2^{n+1}}\ar[d]^-{\text{sign}}\\
\{\pm1\}\ar[r]_-{(-)^2} & \{\pm1\}.
}
\]
The sign~$\rmS_n\to\{\pm1\}$ identifies the abelianisation of the symmetric group, and therefore its first homology group, with the group of order~$2$. This shows that the induced homomorphism in homology is never surjective and, when possible, not injective.
\end{proof}

Those who find the previous example contrived can rest assured: it is not. In fact, the group~$\rmS_{2^n}$ is the full automorphism group of the free Boolean algebra on~$n$ generators, and the stabilisation maps used in the proposition are the canonical ones in that situation. This is explained in detail in~\cite[Sec.~5]{Bohmann--Szymik}, and we shall return to the example in Section~\ref{sec:stable} below.

For those who find permutations and Boolean algebras too nonlinear, we can use a similar argument for the general linear groups over a field~$F$. 

\begin{proposition}\label{prop:GL}
Let 
\[
\GL_{2^0}(F)\longrightarrow\GL_{2^1}(F)\longrightarrow\GL_{2^2}(F)\longrightarrow\dots
\]
be the diagram of general linear groups, where the morphisms are induced by taking an automorphism~$f$ to~$f\otimes\id_{F^2}$. Then the stabilisation homomorphisms
\[
\rmH_1(\GL_{2^n})\to\rmH_1(\GL_{2^{n+1}})
\]
are not injective as soon as~\hbox{$\text{char}(F)\not=2$}, and not surjective as soon as the square class group~$F^\times/F^{\times2}$ is non-trivial.
\end{proposition}

\begin{proof}
From Lemma~\ref{lem:det}, we have~$\det(f\otimes\id_{F^2})=\det(f)^2$. Thus, we have a commutative diagram
\[
\xymatrix@C=4em{
\GL_{2^n}(F)\ar[r]^-{(-)\otimes\id_2}\ar[d]_-{\text{\upshape det}} & \GL_{2^{n+1}}(F)\ar[d]^-{\text{\upshape det}}\\
F^\times\ar[r]_-{(-)^2} & F^\times.
}
\] 
As the determinant induces isomorphisms~$\rmH_1(\GL_n(F))\cong F^\times$ for all~$n\geqslant 1$, the proposition follows from obvious properties of the squaring function: Its kernel is~$\{\pm1\}$, which is non-trivial precisely when~\hbox{$\text{char}(F)\not=2$}, and its cokernel is the square class group~$F^\times/F^{\times2}$, by definition.
\end{proof}
 

\section{Failure of injectivity}\label{sec:inj_fails}

In this section, we explain that the embedding~$\Cr_n(\bbC)\to\Cr_{n+1}(\bbC)$ induces non-injective homomorphisms in homology for large enough values of~$n$. 

\begin{theorem}\label{thm:not_injective}
For any integer~$n\geqslant 3$, the induced homomorphism
\[
\rmH_1(\Cr_n(\bbC);\bbZ)\longrightarrow\rmH_1(\Cr_{n+1}(\bbC);\bbZ)
\]
has a non-trivial kernel.
\end{theorem}

\begin{proof}
Let~$(x,y)=(x,y_1,\dots,y_{n-1})$ be co-ordinates on an~$n$--dimensional rational variety. For every non-zero~$\alpha=\alpha(y)$ in the rational function field~$\bbC(y)$, we can consider the dilation
\[
(x,y)\longmapsto(x\alpha(y),y)
\]
defined by~$\alpha$.

On the one hand, for~$n$ as in the statement of the proposition, Blanc, Lamy, and Zimmermann~\cite[Thm.~A]{BLZ} constructed homomorphisms~$\Cr_n(\bbC)\to\bbZ/2$ and showed that these homomorphisms do not vanish on the elements represented by certain~$\alpha$'s. If we choose such an element~$\alpha$, then the induced homomorphism~\hbox{$\Cr_n(\bbC)^\mathrm{ab}\to\bbZ/2$} from the abelianisation does not vanish on the image of~$\alpha$. Therefore, such an element~$\alpha$ represents a non-zero element in the abelianisation~\hbox{$\Cr_n(\bbC)^\mathrm{ab}\cong\rmH_1(\Cr_n(\bbC);\bbZ)$}.

On the other hand, the stabilisation of the automorphism~$\alpha$ is given by the formula
\begin{equation}\label{eq:stab_of_alpha}
(x,y,z)\longmapsto(x\alpha(y),y,z),
\end{equation}
which leaves~$y$ and the new variable~$z$ fixed. To see that this stabilisation is a commutator, we define
\begin{align*}
c(x,y,z)&=(xz,y,z)\\ 
d(x,y,z)&=(x,y,\alpha(y) z).
\end{align*}
Then, we have~$c^{-1}(x,y,z)=(x/z,y,z)$ and~$d^{-1}(x,y,z)=(x,y,z/\alpha(y))$. We can compute the commutator, therefore, to be
\begin{align*}
[\,c,d\,](x,y,z)
&=cdc^{-1}d^{-1}(x,y,z)\\
&=cdc^{-1}(x,y,z/\alpha(y))\\
&=cd(x\alpha(y)/z,y,z/\alpha(y))\\
&=c(x\alpha(y)/z,y,z)\\
&=(x\alpha(y),y,z).
\end{align*}
These equations show that the stabilisation~\eqref{eq:stab_of_alpha} of the element~$\alpha$ is a commutator in the group of birational equivalences. Thus, it must be zero in the abelianisation~$\rmH_1(\Cr_{n+1}(\bbC);\bbZ)$ of the next Cremona group.
\end{proof}


\section{Failure of surjectivity}\label{sec:sur_fails}

In this section, we explain that the induced map
\[
\rmH_1(\Cr_n(\bbC);\bbZ)\longrightarrow\rmH_1(\Cr_{n+1}(\bbC);\bbZ)
\]
in integral homology is not surjective for large enough values of~$n$. This statement follows from the following result and the universal coefficient theorem, which gives a natural isomorphism~\hbox{$\rmH^1(G;\bbZ)\cong\Hom(\rmH_1(G;\bbZ),\bbZ)$}.

\begin{theorem}\label{thm:not_surjective}
For any integer~$n\geqslant3$, 
the restriction
\[
\rmH^1(\Cr_{n+1}(\bbC);\bbZ)\longrightarrow\rmH^1(\Cr_n(\bbC);\bbZ)
\]
has a non-trivial kernel.
\end{theorem}

The argument here is more involved than that for Theorem~\ref{thm:not_injective}, as it builds on top of deeper results from~\cite{BLZ, BSY} in the context of the Sarkisov program~\cite{Hacon--McKernan}. More recent results announced by Lin and Shinder~\cite{Lin+Shinder:2} strengthen the result to more general bases. Before we prove Theorem~\ref{thm:not_surjective}, we provide some of the geometric background needed to understand the effect of stabilisation.


\subsection{The Sarkisov program}

Faced with the difficulty of describing meaningful generators and relations for groups such as~$\Bir(X)$, Blanc--Lamy--Zimmermann~\cite[p.~270]{BLZ} introduced an equivalent groupoid~$\calB(X)$, i.e., a category in which every morphism is invertible, for which they could provide such a description. We recall the necessary terminology and results from there~(see also~\cite{BSY}).

\begin{definition}
A \textit{Mori fibre space} is a morphism~$X\to B$ with~$\dim(X)>\dim(B)$, relative Picard rank~\hbox{$\rho(X/B)=1$}, and that is the outcome of running a minimal model program: with connected fibres, a normal base, and an anti-canonical divisor~$-K_X$ that is ample over~$B$.
\end{definition}

The objects of the groupoid~$\calB(X)$ are all Mori fibre spaces~$X'\to B$ with~$X'$ birational to~$X$. The morphisms are all birational maps between them. By construction, this groupoid contains~$\Bir(X)$ as the automorphism group of the object~$X$, and, moreover, every object is isomorphic to~$X$, so the groupoid is equivalent to the group when the latter is regarded as a category with one object. It turns out that, in contrast to the group~$\Bir(X)$, the groupoid~$\calB(X)$ does have a meaningful description in terms of generators and relations. To state this description, we need a more general notion that unifies the notions of Mori fibre spaces with the concepts that will describe the generators and relations of the groupoid: Sarkisov links and elementary relations.

\begin{definition} 
A~\textit{rank~$r$ fibration} is a morphism~$\eta\colon X\to B$ that has~\hbox{$\dim(X)>\dim(B)$} and relative Picard rank~\hbox{$\rho(X/B)=r$}. The singularities also have to be reasonable:~$X/B$ is a terminal~$\bbQ$--factorial Mori dream space as defined by Hu--Keel~\cite{Hu--Keel}, the base~$B$ has only Kawamata log terminal singularities, the anti-canonical divisor~$-K_X$ is~$\eta$--big, and for any divisor~$D$ on~$X$, the~$D$--minimal model program from~$X/B$ is terminal~$\bbQ$--factorial.
\end{definition}

Rank~$r$ fibrations are ordered by dominance: we say that one fibration~\hbox{$\eta\colon X\to B$} \textit{dominates} another fibration~\hbox{$\eta'\colon X'\to B'$} when~$\eta$ equals~$\eta'$ composed with a birational contraction~\hbox{$X\to X'$} and a morphism~$B'\to B$ that has connected fibres.


\begin{example}
The rank~$1$ fibrations are the terminal Mori fibre spaces~(see~\cite[Lem.~3.3]{BLZ}). 
\end{example}


\begin{example}
The \hbox{rank~$2$} fibrations correspond to Sarkisov links~(see~\cite[Lem.~3.7]{BLZ}): they factor through exactly two \hbox{rank~$1$} fibrations~$X_j\to B_j$ for~$j=1,2$, that---up to log flips---are divisorial contractions or Mori fibre spaces. This gives four possibilities, depending on whether the morphisms~\hbox{$B_j\to B$} are identities. Each of these cases leads to a diagram
\[
\xymatrix{
X_1\ar[d]\ar@{-->}[rr]^-{\chi}&&X_2\ar[d]\\
B_1\ar[r]&B&B_2,\ar[l]
}
\]
and the birational map~$\chi$ is the {\it Sarkisov link} corresponding to the rank~$2$ fibration. Blanc--Lamy--Zimmermann~\cite[Thm.~4.28~(1)]{BLZ} showed, using the results of Hacon--McKernan~\cite{Hacon--McKernan}, that any morphism in the groupoid~$\calB(X)$ is a composition of automorphisms and Sarkisov links.
\end{example}


\begin{example}
A rank 3 fibration dominates only finitely many Sarkisov links, and these compose to the identity~\cite[Prop.~4.3]{BLZ}. The resulting relations are the \textit{elementary relations}~\cite[Def.~4.4]{BLZ}. Using~\cite{Kaloghiros} and earlier methods developed by Lamy and Zimmermann, they showed that these elementary relations generate all relations~\cite[Thm.~4.28~(2)]{BLZ}.
\end{example}


\subsection{The stabilisation of rank \texorpdfstring{$r$}{r} fibrations}

We will use the following technical lemma in the proof of Theorem~\ref{thm:not_surjective}.

\begin{lemma}\label{lem:lemma}
If~$\eta\colon X\to B$ is a rank~$r$ fibration, then its stabilisation
\[
\eta\times\bbP^1\colon X\times\bbP^1\to B\times\bbP^1
\]
is one, too. In the case~$r=2$, the associated Sarkisov link is~$\chi\times\bbP^1$, where~$\chi$ corresponds to~$\eta$ as above.
\end{lemma}

\begin{proof}[Proof of Lemma~\ref{lem:lemma}]
As~$X$ and~$\bbP^1$ are two~$\bbQ$--factorial varieties, so is their product~\hbox{$X\times\bbP^1$}. If~\hbox{$\xi\colon\widetilde{X}\to X$} is a desingularisation that proves that~$X$ has rational singularities, then we can use the base change~\hbox{$\widetilde{X}\times\bbP^1\to X\times\bbP^1$} to show that~$X\times\bbP^1$ has them, too: the stalks of the higher direct images are the same for~$\xi$ and~$\xi\times\bbP^1$. We can argue similarly for terminal and Kawamata log terminal singularities. These are conditions on the coefficients in the ramification formula
\[
K_{\widetilde{X}}-\pi^*K_X=\sum_ja_jE_j,
\]
and the coefficients are the same for~$\xi\times\bbP^1$:
\begin{align*}
K_{\widetilde{X}\times\bbP^1}-(\xi\times\bbP^1)^*K_{X\times\bbP^1}
&=\pr_1^*(K_{\widetilde{X}}-\xi^*K_X)\phantom{\sum_j}\\
&=\pr_1^*\sum_ja_jE_j\\
&=\sum_ja_j(E_j\times\bbP^1),
\end{align*}
using~$K_{Y\times Z}=\pr_1^*K_Y+\pr_2^*K_Z$ and cancellation of the contributions from~$\bbP^1$ in the difference.

The generic fibres of~$\eta$ and~$\eta\times\bbP^1$ are the same so that any conditions on them~(here: rational connectivity and rational singularities) are inherited from~$\eta$ to~$\eta\times\bbP^1$.
For the condition on the relative Picard number, note that there is a canonical isomorphism~\hbox{$\Pic(X\times\bbP^1)\cong\Pic(X)\oplus\Pic(\bbP^1)$}, and similarly for~$B$, so that the relative Picard numbers stay the same:
\[
\rho(X/B)=\rho(X\times\bbP^1/B\times\bbP^1).
\]
For the second statement, we can apply~`$\times\bbP^1$' to the diagram for~$\chi$ to obtain the diagram for~\hbox{$\chi\times\bbP^1$}.
\end{proof}


\subsection{A proof of Theorem~\ref{thm:not_surjective}}

The main technical result of~\cite{BLZ}, Theorem~D, which is the base for most of the other results in that paper, constructs group homomorphisms~$\Cr_{n+1}(\bbC)\to\bbZ/2$ from Sarkisov links associated with conic bundles. The non-triviality follows from the existence of many different conic bundle models for~$\bbP^n$. This work is extended by Blanc--Schneider--Yasinsky~\cite{BSY}, where the authors exhibit an enormous free quotient of~$\Cr_{n+1}(\bbC)$ and construct group homomorphisms~\hbox{$\Cr_{n+1}(\bbC)\to\bbZ/3$} and~\hbox{$\Cr_{n+1}(\bbC)\to\bbZ$}, this time from Sarkisov~$3$-- and~$6$--links associated with Brauer--Severi fibrations, respectively. 

\begin{proof}[Proof of Theorem~\ref{thm:not_surjective}]
It follows from~\cite[Thm.~B, Thm.~6.2.4, and Sec.~6.4]{BSY} that there exists a Sarkisov~$6$--link~$\chi$ of sufficiently high covering genus that gives rise to a non-zero class in
\[
\rmH^1(\Cr_{n+1}(\bbC);\bbZ)=\Hom(\Cr_{n+1}(\bbC),\bbZ).
\]
To prove the theorem, it now suffices to explain that this class restricts to zero in the previous group~$\rmH^1(\Cr_n(\bbC);\bbZ)=\Hom(\Cr_n(\bbC),\bbZ)$, which means that it evaluates to zero on all elements in $\Cr_{n+1}(\bbC)$ that come from $\Cr_n(\bbC)$.

The class constructed from~$\chi$ has the property that it sends each isomorphism between Mori fibre spaces, each Sarkisov link that is not a~$6$--link, and each~$6$--link that is not equivalent to~$\chi$ or its inverse, to zero~(see~\cite[Thm.~6.2.4]{BSY} again). These facts now imply that the resulting cohomology class will vanish on classes coming from~$\Cr_n(\bbC)$, i.e., on stabilisations of elements in~$\Cr_n(\bbC)$. In fact, the classes coming from~$\Cr_n(\bbC)$ all have~$\bbP^1$--factors in their factorisations by~Lemma~\ref{lem:lemma}: if~$\eta$ is a link with centre~$\Gamma$, then its stabilisation~$\eta\times\bbP^1$ has centre~$\Gamma\times\bbP^1$, which is uniruled, so the covering genus is zero, and it cannot be equivalent to~$\chi$ or its inverse.
\end{proof}


\section{Stable homology}\label{sec:stable}

Unlike in analysis, where convergence is required for a limit to exist, we can consider the stable homology even when homological stability fails. 

\begin{definition}
The \textit{stable homology} of a diagram~\eqref{eq:diagram} of groups is the colimit
\[
\colim\limits_n\rmH_*(G_n;\bbZ).
\]
\end{definition}

In particular, the stable homology of the canonical diagram of the Cremona groups is given as the colimit~$\colim\limits_n\rmH_*(\Cr_n(F);\bbZ)$ of the diagram~\eqref{eq:Cremona_diagram}.

As explained by Quillen~\cite{Quillen:1, Quillen:2}, one reason for our interest in stable homology is that it yields stable characteristic classes, and the rest of this section will present several examples in which it is nontrivial despite the failure of stability. Before we do that, we need to introduce a way to reason about stable homology that we will use in the statements of our results.


\subsection{An infinite loop space}

The hope that the stable invariant is more approachable than the unstable homology is nurtured by standard algebraic K-theory constructions: The path components of a (pointed) space $X$ form a mere (pointed) set $\pi_0(X_0)$. If $X_0=\Omega X_1$ is a loop space, this set is a group, and if~\hbox{$X_0=\Omega X_1=\Omega^2X_2$} is a double loop space, this group is abelian. An infinite loop space requires $X_0=\Omega^n X_n$ for all $n$ in a coherent way. The exposition in~\cite{Adams} is still an excellent introduction to infinite loop spaces and spectra~(see especially~\S3.2). The following result shows that this context applies to the Cremona groups.

\begin{theorem}\label{thm:stable}
For any field, the stable homology of the Cremona groups is the homology of an infinite loop space.
\end{theorem}

\begin{proof}
For a field~$F$, choose a sequence~$x_1,x_2,x_3,\dots$ of variables and let~$\calR_F$ be the category consisting of the~$F$--extensions~$F(n)=F(x_1,\dots,x_n)$ for~\hbox{$n\geqslant0$} and their morphisms~(of~$F$--algebras). This category is a skeleton of the category of all finitely generated rational~$F$--extensions, and we can identify its objects with the natural numbers. 

There is a symmetric monoidal structure~$\boxtimes$ on the category~$\calR_F$ coming from the fact that the~(symmetric monoidal) tensor product of a polynomial ring in~$m$ variables with a polynomial ring in~$n$ variables is a polynomial ring in~$m+n$ variables, and that the field~$F(m+n)$ is a fraction field of the tensor product of the fraction fields~$F(m)$ and~$F(n)$. We write~\hbox{$F(m)\boxtimes F(n)=F(m+n)$} for the induced symmetric monoidal product on the category~$\calR_F$, because it is not the usual tensor product of~$F(m)$ and~$F(n)$ in the category of rings~(for which we have no use here).

This symmetric monoidal structure satisfies cancellation in the sense that the existence of an isomorphism~$F(m_1)\boxtimes F(n)\cong F(m_2)\boxtimes F(n)$ implies~$m_1=m_2$. The unit object~$F(0)=F$ is initial in the category~$\calR_F$ because each~$F(n)$ admits a unique morphism~\hbox{$F\to F(n)$}. Of course, this morphism~$F=F(0)\to F(n)$ can only be an isomorphism if~$n=0$. 

\begin{definition}
The subgroupoid~$\calR_F^\times$ of isomorphisms in the category~$\calR_F$ is the disjoint union of the Cremona groups, thought of as groupoids with one object, and we call~$\calR_F^\times$ the~\textit{Cremona groupoid}. 
\end{definition}

We can now apply Quillen's algebraic K-theory construction~\cite{Grayson} that associates to any such symmetric monoidal groupoid a group completion, which is an infinite loop space whose homology computes the stable homology of the automorphism groups. We refer to~\cite[Sec.~2]{Bohmann--Szymik} and~\cite[Sec.~3.2]{RWW} for details. For the Cremona groupoid~$\calR_F^\times$, we get the infinite loop space we were looking for.
\end{proof}


Let us write~$\rmK(\calR^\times_F)$ for the K-theory spectrum associated with the Cremona groupoid. The underlying infinite loop space is denoted by~$\Omega^\infty\rmK(\calR^\times_F)$. Because of cancellation, its abelian group of components is given by
\[
\rmK_0(\calR^\times_F)=\pi_0\Omega^\infty\rmK(\calR^\times_F)=\bbZ.
\]
The component~$\Omega^\infty_0\rmK(\calR^\times_F)$ of the zero element describes the stable homology of the Cremona groups:
\[
\colim\limits_n\rmH_*(\Cr_n(F);\bbZ)=\rmH_*\Omega^\infty_0\rmK(\calR^\times_F).
\]
We refer to~\cite[Sec.~2]{Bohmann--Szymik} and~\cite[Sec.~3.2]{RWW} again.

\begin{remark}
If~$F\to E$ is a field extension, we have a morphism
\[
\rmK(\calR^\times_F)\longrightarrow\rmK(\calR^\times_E),
\]
and we can take that to indicate that the cases where~$F$ is a prime field are the most fundamental ones.
\end{remark}

We can now turn to the question of how to describe the infinite loop space of Theorem~\ref{thm:stable}, or at least show that it is non-trivial, i.e., not a $\bbZ$'s worth of contractible components. Before we do that, we build intuition by revisiting the examples from Section~\ref{sec:homological_instability} from the perspective of stable homology. 


\subsection{Revisiting the examples from Section~\ref{sec:homological_instability}}\label{sec:ex2rev}

The following two remarks show that the stable homology is non-trivial in the two examples of homological instability discussed in Section~\ref{sec:homological_instability}.

\begin{remark}
If we consider the symmetric monoidal groupoid of finite sets and bijections with respect to disjoint union, it is clear that the components give the additive monoid of natural numbers, and stabilisation would happen with respect to the canonical generator. In contrast, when we consider the symmetric monoidal groupoid of non-empty finite sets and bijections with respect to the cartesian product, then the abelian monoid of components is given by the positive integers with respect to multiplication, and this is generated~(freely) by the primes, and therefore, has no preferred generator. Consequently, there are several choices for stabilisation. For Proposition~\ref{prop:Boole}, we stabilised the product with the element representing~$2$. Then the proof of Proposition~\ref{prop:Boole} shows that the first stable homology is
\[
\colim(\bbZ/2\overset{\cdot2}{\longrightarrow}\bbZ/2\overset{\cdot2}{\longrightarrow}\bbZ/2\overset{\cdot2}{\longrightarrow}\dots)=(\bbZ/2)\left[\frac12\right]=0.
\]
Still, even though the~$2$--torsion disappears, there is non-trivial higher stable homology in this case. To see this, we can identify the diagram in Proposition~\ref{prop:Boole} with the diagram
\[
\Aut(B_0)\longrightarrow\Aut(B_1)\longrightarrow\Aut(B_2)\longrightarrow\dots,
\]
where~$B_n$ is the free Boolean algebra on~$n$ generators. In that language, the stable homology has been computed in~\cite[Sec.~5]{Bohmann--Szymik}. It is given by~$\Omega^\infty_0\bbS[1/2]^\times$, where~$\bbS$ is the sphere spectrum, whose homotopy groups are the stable homotopy groups of spheres. The~$2$--torsion has been localised away, and the first non-zero homotopy group is~$\pi_3$, which has order~$3$ and is generated by the stable Hopf map~$\nu$.
\end{remark}


\begin{remark}\label{rem:GL-2}
We can make a similar statement, mutatis mutandis, for vector spaces. Looking back at Proposition~\ref{prop:GL}, we see that the stable homology computation for the diagram of the general linear groups over a field~$F$ with respect to~$\otimes F^2$--stabilisation is similar. The process inverts multiplication with the class~$2=[F^2]$ in the K-theory spectrum~$\rmK(F)$, and the stable homology is given as in~\cite{Bohmann--Szymik} by the components of the infinite loop space~$\Omega^\infty_0\rmK(F)[1/2]^\times$ that correspond to the units in the ring~\hbox{$\pi_0\rmK(F)[1/2]^\times=\bbZ[1/2]^\times$} that are powers of~$2$. The more drastic case, where \textit{all} the classes~$n=[F^n]$ for~$n\geqslant2$ are inverted, goes back to~\cite{May}~(see~\cite[Prop.~5]{Weibel}), and involves---now unsurprisingly---a rationalisation. Even in that case, the higher stable homology can still be non-trivial, as the K-theory of fields is often non-trivial, even rationally.
\end{remark}


\subsection{Producing stable homology classes}

After the digression of Section~\ref{sec:ex2rev}, we now return to~$\rmK(\calR^\times_F)$ and first describe various ways to map into it.

\subsubsection{K-theory classes from automorphisms}

The automorphism group of any object~$X$ in a symmetric monoidal groupoid~$(\calG,\boxtimes)$ gives rise to compatible morphisms
\[
\Aut(X)\wr\rmS_n\longrightarrow\Aut(X^{\boxtimes n})
\]
of groups that combine to give a morphism
\[
\Sigma^\infty_+\rmB\Aut(X)\longrightarrow\rmK(\calG)
\]
of spectra. 

In particular, for~$\calG=\calR^\times_F$ and any rational~$F$--variety~$X$, this method gives rise to a morphism
\[
\Sigma^\infty_+\rmB\Bir(X)\longrightarrow\rmK(\calR^\times_F), 
\]
and precomposing with the morphism induced by the inclusion $\Aut(X)\to\Bir(X)$ of groups, we get 
\[
\Sigma^\infty_+\rmB\Aut(X)\longrightarrow\rmK(\calR^\times_F),
\]
which can be advantageous when $\Aut(X)$ is easier to handle than~$\Bir(X)$. Note that both $\Aut(X)$ and~$\Bir(X)$ are automorphism groups of $X$, albeit in different categories!

In particular, if we specialise to the projective line~$X=\bbP^1_F$, where every birational equivalence is an automorphism, we get 
\begin{equation}
\Sigma^\infty_+\rmB\PGL_2(F)\longrightarrow\rmK(\calR^\times_F).
\end{equation}


\subsubsection{Affine spaces and the additive group}

The automorphism groups of affine space~$\bbA^n_F$ are big and complicated. However, as~$\bbA^n_F$ is the underlying space of the group~$(\GGa^n)_F$, we can focus on the automorphisms that respect the group structure, which are given by the general linear groups~$\GL_n(F)\leqslant\Aut(\bbA^n_F)$. The inclusions are compatible with the canonical stabilisation maps, and this leads to a morphism
\begin{equation}\label{eq:K(F)}
\rmK(F)\longrightarrow\rmK(\calR^\times_F)
\end{equation}
from the algebraic K-theory spectrum~$\rmK(F)$ of the field~$F$ into our~K-theory spectrum. One should think of this as a highly conceptual way of saying that the Cremona groups contain the linear changes of coordinates.


\subsubsection{Tori and the multiplicative group}

We can proceed similarly for the tori~$(\GGm^n)_F$, which are also rational and carry a group structure. Here, we can note that~\hbox{$\GGm^n=\Hom(\bbZ^n,\GGm)$}, so that the groups~$\GL_n(\bbZ)$ act on them.
The homomorphisms~\hbox{$\GL_n(\bbZ)\to\Cr_n(F)$} are compatible with each other and lead to a morphism
\begin{equation}\label{eq:K(ZZ)}
\rmK(\bbZ)\longrightarrow\rmK(\calR^\times_F)
\end{equation}
of spectra. One should think of this as a highly conceptual way of exploiting Demazure's idea~\cite[p.~522]{Demazure} that the~$\GL_n(\bbZ)$ are the Weyl groups of the maximal tori that sit inside the Cremona groups as groups of monomial transformations.

\begin{remark}
Even though there is a canonical morphism~$\rmK(\bbZ)\to\rmK(F)$ from the unique ring map~\hbox{$\bbZ\to F$}, there is no obvious reason why the restriction of~\eqref{eq:K(F)} along this morphism should give~\eqref{eq:K(ZZ)}. Just looking at the definitions, the corresponding subgroups of the Cremona groups are very different: one is made of monomial transformations, the other of linear transformations.
\end{remark}


\subsection{Detecting stable homology classes}

So far, we have described many ways to map into~$\rmK(\calR^\times_F)$. While this may yield many homotopy classes coming from the various sources, it does not prove that the target~$\rmK(\calR^\times_F)$ is non-trivial, where trivial means discrete, so it has no higher homotopy groups. The difficulty lies in finding systematic morphisms from the Cremona groups into other groups. The ones we used for homological instability are, by their very nature, not suitable for this purpose. Instead, we rely on~\cite{Lin+Shinder:1}~(see also~\cite{GLU}) to find the following result:

\begin{theorem}\label{thm:non-trivial}
The first homotopy group~\hbox{$\rmK_1(\calR^\times_F)=\pi_1\rmK(\calR^\times_F)$} is rationally non-zero for any infinite field~$F$. In particular, the spectrum~$\rmK(\calR^\times_F)$ is not discrete 
\end{theorem}

\begin{proof}
We can define~$\Burn_n(F)$ to be the free abelian group on the set of isomorphism classes of field extensions of~$F$ that have transcendental degree~$n$, or as in~\cite{Kontsevich+Tschinkel}, with birational classes of~$F$--varieties instead, whichever we find more convenient. By passing from a variety to its function field, we see that both constructions are isomorphic. 

For every~$n$--dimensional~$F$--variety~$X$, Lin and Shinder~\cite{Lin+Shinder:1} constructed a homomorphism
\[
c\colon\Bir(X)\longrightarrow\Burn_{n-1}(F)
\]
by sending a birational equivalence to the difference of~(the function fields of) its exceptional divisors. By~\cite[Thm.~4.4]{Lin+Shinder:1}, there are two non-isomorphic Calabi--Yau~$3$--folds~$Z$ and~$Z'$, and a birational equivalence~$f\in\Bir(\bbP^5_F)$ such that
\[
c(f)=[Z\times_F\bbP^1_F]-[Z'\times_F\bbP^1_F]\not=0.
\]

It follows inductively from the commutativity of the diagram
\[
\xymatrix{
\Bir(X)\ar[d]_-{f\longmapsto f\times\id_{\bbP^1_F}}\ar[r]^-c &\Burn_{n-1}(F)\ar[d]^-{[E]\longmapsto[E
\times_F\bbP^1_F]}\\
\Bir(X\times_F\bbP^1_F)\ar[r]_-c &\Burn_n(F),
}
\]
implicit in the proof of their result, that
\begin{equation}\label{eq:stable_char_class}
c(f\times\id_{\bbP^1_F}\times\dots\times\id_{\bbP^1_F})=c(f)\times_F\bbP^1_F\times_F\dots\times_F\bbP^1_F.
\end{equation}
Since the Calabi--Yau~$3$--folds~$Z$ and~$Z'$ are not stably birational by~\cite[Cor.~A.5]{Lin+Shinder:1}, the left-hand side is nonzero regardless of the number of factors~$\bbP^1_F$.

If we define homomorphisms~$\Cr_n(F)\to\bbZ$ by sending~$f$ to the difference of the coefficients of~$c(f)$ in front of~$[F(Z\times_F(\bbP^1_F)^{n-3})]$ and~\hbox{$[F(Z'\times_F(\bbP^1_F)^{n-3})]$}, these are compatible by~\eqref{eq:stable_char_class} and define a stable characteristic class in degree~$1$, and this class is nonzero by the properties of the birational equivalence~$f$ and the Calabi--Yau~$3$--folds~$Z$ and $Z'$ mentioned before. This shows~$\lim_n\rmH^1(\Cr_n(F);\bbZ)\not=0$.

To drive the result home, we deduce from
\[
\Hom(\colim_n\rmH_1(\Cr_n(F)),\bbZ)=\lim_n\rmH^1(\Cr_n(F);\bbZ)\not=0,
\]
that
\[
\rmH_1\Omega^\infty_0\rmK(\calR^\times_F)=\rmH_1(\colim_n\rmB\Cr_n(F))=\colim_n\rmH_1(\rmB\Cr_n(F))
\]
must be non-zero. By Hurewicz, this agrees with $\pi_1\Omega^\infty_0\rmK(\calR^\times_F)=
\rmH_1\Omega^\infty_0\rmK(\calR^\times_F)$, and the homotopy groups of an infinite loop space agree with those of the corresponding spectrum in positive dimensions.
\end{proof}

The homotopy groups~$\rmK_i(\calR^\times_F)$ are invariants of the field~$F$. It would certainly be a challenging endeavour to express these in terms of other invariants of the field~$F$.


%
%
%


%



\begin{thebibliography}{99}

\bibitem{Adams} J.F. Adams. Infinite loop spaces. Annals of Mathematics Studies, 90. Princeton University Press, 1978.

\bibitem{Arnold} V.I. Arnold. On some topological invariants of algebraic functions. Tr. Mosk. Mat. Obs. 21 (1970) 27--46.

\bibitem{Bohmann--Szymik} A.M. Bohmann, M. Szymik. Boolean algebras, Morita invariance and the algebraic K-theory of Lawvere theories. Math. Proc. Cambridge Philos. Soc. 175 (2023) 253--270.

\bibitem{BF13} J. Blanc, J.-P. Furter. Topologies and structures of the Cremona groups. Ann. of Math. 178 (2013) 1173--1198. 

\bibitem{BLZ} J. Blanc, S. Lamy, S. Zimmermann. Quotients of higher-dimensional Cremona groups. Acta Math. 226 (2021) 211--318.

\bibitem{BSY} J. Blanc, J. Schneider, E. Yasinsky. Birational maps of Severi--Brauer surfaces, with applications to Cremona groups of higher rank. Duke Math.~J.~(to appear)

\bibitem{BZ18} J. Blanc, S. Zimmermann. Topological simplicity of the Cremona groups. Amer. J. Math. 140 (2018) 1297--1309.

\bibitem{Demazure} M. Demazure. Sous-groupes alg\'ebriques de rang maximum du groupe de Cremona. Ann. Sci. \'Ecole Norm. Sup. 3 (1970) 507--588. 

\bibitem{EVW} J.S.~Ellenberg, A.~Venkatesh, C.~Westerland. Homological stability for Hurwitz spaces and the Cohen--Lenstra conjecture over function fields. Ann. of Math. 183 (2016) 729--786.

\bibitem{GLU} A. Genevois, A. Lonjou, C. Urech. On a theorem by Lin and Shinder through the lens of median geometry. Selecta Math. 31 (2025) No.~1, Paper 17, 10 pp.

\bibitem{Grayson} D. Grayson. Higher algebraic K-theory. II (after Daniel Quillen). Algebraic K-theory (Proc. Conf., Northwestern Univ., Evanston, Ill., 1976) 217--240. Lecture Notes in Math. 551. Springer-Verlag, Berlin-New York, 1976.

\bibitem{Hacon--McKernan} C.D. Hacon, J. McKernan. The Sarkisov program. J. Algebraic Geom. 22 (2013) 389--405.

\bibitem{Harer} J.L. Harer. Stability of the homology of the mapping class groups of orientable surfaces. Ann. of Math. 121 (1985) 215--249.

\bibitem{Hatcher} A. Hatcher. Homological stability for automorphism groups of free groups. Comment. Math. Helv. 70 (1995) 39--62.

\bibitem{Hu--Keel} Y. Hu, S. Keel. Mori dream spaces and GIT. Dedicated to William Fulton on the occasion of his 60th birthday. Michigan Math. J. 48 (2000) 331--348.

\bibitem{Kaloghiros} A.-S. Kaloghiros. Relations in the Sarkisov program. Compos. Math. 149 (2013) 1685--1709.

\bibitem{Kontsevich+Tschinkel} M. Kontsevich, Y. Tschinkel. Specialization of birational types. Invent. Math. 217 (2019) 415--432.

\bibitem{Lin+Shinder:1} H.-Y. Lin, E. Shinder. Motivic invariants of birational maps. Ann. of Math. 199 (2024) 445--478.

\bibitem{Lin+Shinder:2} H.-Y. Lin, E. Shinder. Unboundedness for motivic invariant of birational automorphisms. \href{https://arxiv.org/abs/2510.00290}{arXiv:2510.00290} 

\bibitem{May} J.P. May, F. Quinn, N. Ray, and J. Tornehave. E$_\infty$ ring spaces and E$_\infty$ ring spectra. Lecture Notes in Mathematics, 577. Springer-Verlag, Berlin-New York, 1977. 

\bibitem{Nakaoka} M. Nakaoka. Decomposition theorem for homology groups of symmetric groups. Ann. of Math. 71 (1960) 16--42.

\bibitem{Quillen:1} D. Quillen. Notebooks, 1971. 
\href{https://www.claymath.org/online-resources/quillen-notebooks}{www.claymath.org/online-resources/quillen-note\-books}

\bibitem{Quillen:2} D. Quillen. Characteristic classes of representations. Algebraic K-theory (Proc. Conf., Northwestern Univ., Evanston, Ill., 1976) 189--216. Lecture Notes in Math., Vol. 551. Springer, Berlin-New York, 1976. 

\bibitem{RWW} O. Randal-Williams, N. Wahl. Homological stability for automorphism groups. Adv. Math. 318 (2017) 534--626.

\bibitem{Szymik} M. Szymik. Twisted homological stability for extensions and automorphism groups of free nilpotent groups. J. K-Theory 14 (2014) 185--201. 

\bibitem{Szymik+Wahl} M. Szymik, N. Wahl. The homology of the Higman--Thompson groups. Invent. Math. 216 (2019) 445--518.

\bibitem{Wahl} N. Wahl. Homological stability: a tool for computations. Proc. Int. Cong. Math. 2022, Vol. 4, 2904--2927. Ed. D. Beliaev and S. Smirnov. European Mathematical Society, 2023.

\bibitem{Weibel} C.A. Weibel, K-theory of Azumaya algebras. Proc. Amer. Math. Soc. 81 (1981) 1--7. 

\end{thebibliography}
\end{document}